\documentclass[12pt]{amsart}
\usepackage[margin=1.5in]{geometry}
\usepackage{amsthm,amsmath,amssymb}
\theoremstyle{plain}
\newtheorem{theorem}{Theorem}[section]
\newtheorem{lemma}[theorem]{Lemma}
\newtheorem{question}[theorem]{Question}
\numberwithin{equation}{section}
\newcommand{\gen}[1]{\langle#1\rangle}
\newcommand{\abs}[1]{\lvert#1\rvert}
\newcommand*{\set}[1]{\{#1\}}

\newcommand{\FF}{\mathbb F}
\DeclareMathOperator{\GL}{GL}

\DeclareMathOperator{\PSL}{PSL}

\DeclareMathOperator{\soc}{soc}

\usepackage{xcolor} % purple
\usepackage{url}
\usepackage[T1]{fontenc}
\usepackage{lmodern}
\usepackage{listings}
\definecolor{kw}{RGB}{0,0,150}
\definecolor{str}{RGB}{150,0,0}
\definecolor{cmt}{RGB}{60,120,60}

\lstdefinelanguage{Sage}{
  language=Python,
  morekeywords={matrix,GF,identity,diagonal}
}
\usepackage[
  backend=biber,
  style=numeric,
  doi=true,
  url=false,
  giveninits=true,
  isbn=false
]{biblatex}
\usepackage[colorlinks,linkcolor=blue,citecolor=purple]{hyperref}
\title{On transitive sets of derangements in primitive groups}
\author{Peter M\"uller}
\address{Institute of Mathematics, University of W\"urzburg}
\email{peter.mueller@uni-wuerzburg.de}
\begin{document}
\begin{abstract}
  We construct a primitive permutation action of the Steinberg
  triality group $^3D_4(2)$ of degree $4{,}064{,}256$ and show that
  there are distinct points $\alpha,\beta$ such that there is no
  derangement $g\in{{}^3D_4}(2)$ with $\alpha^g=\beta$. This answers a
  question by John G.~Thompson (Problem 8.75 in the Kourovka notebook
  \cite{kourovka_notebook}) in the negative.
\end{abstract}
\maketitle
\section{Introduction}
Beginning with the 8th edition in 1982, the Kourovka Notebook
\cite{kourovka_notebook}, which collects unsolved problems in group
theory, has included the following question -- attributed to John
G.~Thompson and designated ``A known problem'' -- as Problem 8.75:
\begin{question}\label{q:8.75}
  Suppose $G$ is a finite primitive permutation group on $\Omega$, and
  $\alpha$, $\beta$ are distinct points of $\Omega$. Does there exist
  an element $g\in G$ such that $\alpha^g=\beta$ and $g$ fixes no
  point of $\Omega$?
\end{question}
A potential interest of this is that if the answer were affirmative
with a proof not using character theory, then Frobenius' theorem that
the kernel of a finite Frobenius group is a subgroup would admit a
proof devoid of character theory. Such a proof still does not exist,
despite some progress like \cite{flavell___frobenius} and
\cite{tao___frobenius}. See the argument in \cite[Assertion
2]{kaplan_lev___transversals} (which relies on
\cite{flavell___frobenius}) how Frobenius' theorem would follow from a
positive answer of Question \ref{q:8.75}.

We refer to the survey \cite{burness_giudici___classical} for more
open problems (like this one in
\cite[Remark~1.1.3(iii)]{burness_giudici___classical}) and results
about derangements in permutation groups.

For most Aschbacher-ONan-Scott types of primitive groups the answer to
Question \ref{q:8.75} is trivially yes in view of the presence of
regular subgroups. As we considered it unlikely to find a
counterexample for primitive groups of product type (Class III in
\cite{liebeck_praeger_saxl___onan}), our search focused on the almost
simple action of groups and eventually produced a counterexample:
\begin{theorem}\label{t:3d4}
  The Steinberg triality group $G = {}^{3}D_{4}(2)$ of order
  $211{,}341{,}312 = 2^{12}\cdot 3^{4}\cdot 7^{2}\cdot 13$ admits a
  primitive permutation action on a set $\Omega$ of size
  $4{,}064{,}256$ with the following property: There exist distinct
  elements $\alpha, \beta\in\Omega$ such that every $g\in G$ with
  $\alpha^g=\beta$ has at least one fixed point in $\Omega$.
\end{theorem}
This counterexample is the only one we know of. Let $\soc(G)$ be the
socle of an almost simple group $G$. There are $584$ nonabelian simple
groups $S$ of order less than or equal $20{,}158{,}709{,}760$. A brute
force computation using Magma \cite{magma}, where $S=\soc(G)$ is one
of these simple groups, resulted in the only example from Theorem
\ref{t:3d4} . (The cases $S=\PSL(2,q)$ were only checked for prime
powers $q\le 1789$. But this case can be ruled out anyway, see Remark
(g).) We also used the data base of primitive permutation groups in
Magma \cite{magma} to check that there is no counterexample of degree
less than $8{,}192$.

In contrast to the demanding search for a counterexample, the proof of
Theorem \ref{t:3d4} will require little computation. We will work with
${}^{3}D_{4}(2)$ as a subgroup of $\GL_8(\FF_8)$. Except for some
general group theoretic arguments and the known fact that
${}^{3}D_{4}(2)$ has a maximal subgroup isomorphic to
$C_{13}\rtimes{C_4}$, the verification of Theorem \ref{t:3d4} is
reduced to check some identities for specific matrices from
$\GL_8(\FF_8)$.
\section{Proof of Theorem \ref{t:3d4}}
The proof will be based on the following lemma.
\begin{lemma}\label{l:mat8}
  The Steinberg triality group $G = {}^{3}D_{4}(2)$ contains elements
  $s, x$, and $y$ with the following properties:
  \begin{itemize}
  \item[(1)] $x$, $y$, and $yx$ are conjugate in $G$,
  \item[(2)] $x^y = x^{-1}$,
  \item[(3)] $x^2=y^2$,
  \item[(4)] $s$ has order $13$,
  \item[(5)] $s^y = s^8$,\label{i:e}
  \item[(6)] $x\notin\gen{s, y}$,
  \item[(7)] $(hx)^{13}=1$ for every $h\in\gen{s, y}\setminus\gen{y}$.
  %\item[()] $(hx)^2 = y^2$ for every $h\in\gen{y}$.
  \end{itemize}
\end{lemma}
Before proving this lemma, we show how Theorem \ref{t:3d4} follows from
it. By (4) and (5), $y$ acts as an automorphism of order $4$ on the
Sylow $13$-subgroup $S=\gen{s}$. According to \cite{Kleidman:3D4}, $G$
has a maximal subgroup isomorphic to $C_{13}\rtimes C_{4}$. Thus $y$
has order $4$, and $G$ acts primitively on the
$\abs{G}/52=4{,}064{,}256$ right cosets of
$H=\gen{s, y}=\gen{s}\rtimes\gen{y}$ in $G$.

Note that $H\ne Hx$ by (6). We claim that every $g\in Hx$ fixes a
right coset of $H$. Write $g=hx$ for some $h\in H$. We need to show
that $hx$ is contained in some conjugate $H^z$ of $H$, for then
$Hzg=Hzhx=Hz$.

If $h\in\gen{s, y}\setminus\gen{y}$, then $hx$ lies in a conjugate of
$S<H$ by (7).

It remains to look at the case $h\in\gen{y}$, that is
$hx\in\set{x, yx, y^2x, y^{-1}x}$. By (1), $x$ and $yx$ are conjugate
to $y\in H$. From (2) we get $y^{-1}x=x^{-1}y^{-1}=(yx)^{-1}$, so
$y^{-1}x$ is conjugate to $y^{-1}\in H$ by (1). Finally, $y^2x=x^3$ by
(3), so $y^2x$ is conjugate to $y^3\in H$ by (1).
\subsection{Proof of Lemma \ref{l:mat8}.}
We work with a representation of $G={}^3D_4(2)$ as a subgroup of
$\GL_8(\FF_8)$ as given in \cite[Section
4.2]{howlett_rylands_taylor}: Pick $\omega\in\FF_8$ with
$\omega^3+\omega+1=0$. Note that $\omega$ has multiplicative order
$7$. Let $E_{i,j}$ be the $8\times 8$ matrix over $\FF_8$ with $1$ in
the $i$th row and $j$th column and $0$ elsewhere. Furthermore, let $E$
be the $8\times 8$ identity matrix. Using the notation from
\cite{howlett_rylands_taylor}, we set
\begin{align*}
  x_R(1) &=  E + E_{1,2} + E_{3,4} + E_{3,5} + E_{3,6} + E_{4,6} +%
           E_{5,6} + E_{7,8}\\
  n &= E_{1,3} + E_{2,1} + E_{3,7} + E_{4,5} + E_{5,4} +%
      E_{6,2} + E_{7,8} + E_{8,6}\\
  h_R(\omega) &= \text{diagonal matrix with diagonal }(\omega^4, \omega^3, \omega^3,%
                \omega, \omega^6, \omega^4, \omega^4, \omega^3).
\end{align*}
Then $G=\gen{a, b}$, where $a=x_R(1)\cdot n$ and $b=h_R(\omega)$.

Essentially the same generators $a$ and $b$ are given in
\cite[\href{https://www.lmfdb.org/Groups/Abstract/211341312.a}{Abstract
  group 211341312.a}]{lmfdb}, the only difference is that the fourth
and fifth rows and columns are switched there.

With these matrices $a$ and $b$ set
\begin{align*}%
  x &= b\cdot (a^{2}\cdot b^{2}\cdot a)^{2}\cdot b^{-2}\cdot a\cdot
      b^{-2}\\
  t &= a^{-3}\cdot b\cdot (b\cdot a)^{2}\cdot b^{-2}\cdot a^{2}\cdot
      b\\
  y &= x^t\\
  s &= a\cdot a^b\\
  d &= x\cdot s\cdot y\cdot s\cdot a\cdot y\cdot s^{2}\cdot t
\end{align*}      
With these settings, the following SageMath \cite{sagemath} code
verifies all the claims instantly (in less than a second). The
comments indicate the respective items from Lemma \ref{l:mat8}. If one
wants to confirm that the matrices $a$ and $b$ are correctly copied
(and correctly given in \cite{howlett_rylands_taylor}), one may
uncomment the last line. Note that except for the last line, the code
merely does simple matrix operations which in principle could be done
by hand.%
\begin{lstlisting}
F.<w> = GF(8, modulus=[1, 1, 0, 1])
e = matrix.identity(F, 8)
xr1_ij = [(1, 2), (3, 4), (3, 5), (3, 6), (4, 6), (5, 6), (7, 8)]
n_ij = [(1, 3), (2, 1), (3, 7), (4, 5), (5, 4), (6, 2), (7, 8), (8, 6)]
xr1 = e + matrix(F, 8, {(i-1, j-1): 1 for i, j in xr1_ij})
n = matrix(F, 8, {(i-1, j-1): 1 for i, j in n_ij})
a = xr1 * n
b = matrix.diagonal(F, [w^z for z in [4, 3, 3, 1, 6, 4, 4, 3]])
x = b * (a^2 * b^2 * a)^2 * b^-2 * a * b^-2
t = a^-3 * b * (b * a)^2 * b^-2 * a^2 * b
y = t^-1 * x * t                                 # one part of (1)
s = a * b^-1 * a * b
d = x * s * y * s * a * y * s^2 * t
S = [s^i for i in range(13)]
Y = [y^i for i in range(4)]
H = [u * v for u in S for v in Y]
assert d^-1 * x * d == y * x                     # other part of (1)
assert y^-1 * x * y == x^-1                      # (2)
assert x^2 == y^2                                # (3)
assert s != e and s^13 == e                      # (4)
assert y^-1 * s * y == s^8                       # (5)
assert x not in H                                # (6)
assert all(h in Y or (h * x)^13 == e for h in H) # (7)
print('All claims verified')
# Uncomment to optionally check that a and b indeed generate 3D(4,2):
# assert MatrixGroup([a, b]).structure_description() == '3D(4,2)'
\end{lstlisting}
This SageMath code is provided at \cite{dera_3d42_verify}, where it
can be downloaded and run online.
\section{Remarks}
\begin{itemize}
\item[(a)] If we replace $G$ by its automorphism group ${}^3D_4(2).3$,
  then the primitive action of degree $4{,}064{,}256$ has a transitive
  set of derangements, so it is not a counterexample to Question
  \ref{q:8.75}.
\item[(b)] Let $H$ be a subgroup of a finite group $G$. Then the
  inclusion $Hx\subseteq\cup_{g\in G}H^g$ holds for $x\in G$ if and
  only if $HxH\subseteq\cup_{g\in G}H^g$. For $G$ and a point
  stabilizer $H$ as in Theorem \ref{t:3d4}, there are precisely
  $78{,}366$ double cosets $HxH$ different from $H$, and it turns out
  that exactly one of these double cosets accounts for the
  counterexample!
\item[(c)] In hindsight, this kind of uniqueness of the counterexample
  helped to derive some of the structural properties in Lemma
  \ref{l:mat8}. For instance, if $HxH$ gives a counterexample, then so
  does $Hx^{-1}H$, thus $HxH=Hx^{-1}H$ by uniqueness. So
  $x^{-1}=h_1xh_2$ for $h_i\in H$. This gives
  $(xh_1)^{-1}=(xh_1)(h_1^{-1}h_2)$. So upon replacing $x$ with
  $xh_1$, we may assume $x^2\in H$. Now $x$ has order $\ne13$, for
  otherwise $x=(x^2)^7\in H$. Thus $x$ has order $2$ or $4$. One can
  rule out the former case. Again, let $S$ be the Sylow $13$-subgroup
  of $H$. The normalizer of $S$ in the automorphism group
  ${}^3D_4(2).3$ of $G$ has the form $S\rtimes\gen{\sigma}$ with
  $\sigma$ of order $12$. Clearly, $\sigma$ normalizes $H$, and we may
  assume that $\sigma^6=x^2$. Here $\sigma^3$ plays the role of $y$ in
  Lemma \ref{l:mat8}. From $Hx^\sigma H=HxH$ we get further properties
  from Lemma \ref{l:mat8}.
\item[(d)] I have tried to extend the counterexample to the bigger
  twisted triality groups $G={}^3D_4(q)$ for $q=3$ and $q=4$ and the
  maximal subgroup $H=C_{q^4-q^2+1}\rtimes C_4$. Here, $G$ is far too
  large for it to be feasible to compute representatives $x$ of the
  double cosets of $H$ in $G$. I was only able to verify that the
  elements $x$ of order $4$ with $x^2\in H$ do not produce
  counterexamples (unlike in the case $q=2$).
\item[(e)] Initially, I had tried to explore a vast generalization of
  Question \ref{q:8.75}. The primitivity of a group can not be read
  off from the permutation character, see the negative answer
  \cite{guralnick_saxl___primitive_permutation_characters} by
  Guralnick and Saxl to a question of Wielandt. However, by a result
  of Higman (see \cite{higman___intersection} or \cite[Theorem
  3.2A]{DixMort}), primitivity is encoded in the orbital graphs: The
  group $G$ is primitive if and only if each non-diagonal orbital
  graph is connected.

  Now there was some experimental evidence that the following
  generalization of Question \ref{q:8.75} could hold: Let $G$ act
  transitively on the finite set $\Omega$. Pick distinct
  $\alpha,\beta$ in $\Omega$. If the orbital graph which contains the
  edge $(\alpha, \beta)$ is connected, then there is a derangement
  $g\in G$ with $\alpha^g=\beta$. The smallest counterexample which I
  found is $G=\PSL_3(4)$ acting on $4032$ points (so the point
  stabilizer has order $5$).
\item[(f)] I have some evidence (and partial results) that the
  following could hold: Let $G$ act transitively on the finite set
  $\Omega$, and pick distinct $\alpha,\beta$ in $\Omega$. Then there
  is an element $g\in G$ with $\alpha^g=\beta$ whose number of fixed
  points is different from $1$. This is now Problem 21.99 in the 21st
  edition (2026) of the Kourovka Notebook \cite{kourovka_notebook}.

  Applied to a finite Frobenius
  permutation group acting on a finite set $\Omega$, this would mean
  that for any distinct $\alpha, \beta\in\Omega$, there is a unique
  derangement $g\in G$ with $\alpha^g=\beta$. But then Frobenius'
  theorem follows immediately from a simple counting argument (see
  e.g.~\cite[Assertion 1]{kaplan_lev___transversals}).
\item[(g)] It would still be interesting to decide for which primitive
  groups Question \ref{q:8.75} has a positive answer. With some
  effort, one can show that this is the case if $G$ is almost simple
  with socle $\soc(G)=\PSL_2(\FF_q)$. See \cite{mueller___psl2q}.
\end{itemize}
\printbibliography
\end{document}